\documentclass[11pt,a4paper]{article}
\usepackage[affil-it]{authblk}
\usepackage{graphicx} 
\usepackage[margin=1in]{geometry}
\usepackage{amsmath}
\usepackage{amsthm}
\usepackage{caption}
\usepackage{algorithm}
\usepackage{algpseudocode}
\usepackage{amsfonts}
\usepackage{xcolor}
\usepackage{array}

\newtheorem{theorem}{Theorem}
\newtheorem{remark}[theorem]{Remark}

\newtheorem{prop}[theorem]{Proposition}

\newtheorem{definition}{Definition}
\newtheorem{examp}{Example}
    	  \newcommand{\bA}{\mbox{\boldmath${A}$}}
            \newcommand{\bJ}{\mbox{\boldmath${J}$}}

\newcommand{\bB}{{\bf B}}

\newcommand{\bC}{{\bf C}}

\newlength{\figheight}
\newcommand{\bI}{{\bf I}}

\newcommand{\bQ}{{\bf Q}}
\newcommand{\bR}{{\bf R}}
\newcommand{\bS}{{\bf S}}

            \newcommand{\bV}{\mbox{\boldmath${V}$}}

\newcommand{\be}{{\bf e}}

\newcommand{\bx}{{\bf x}}

\newcommand{\hu}{\hat{u}}
\newcommand{\hx}{\hat{x}}
\newcommand{\hy}{\hat{y}}

\newcommand{\RR}{{\mathbb{R}}}
\newcommand{\CC}{{\mathbb{C}}}
\def\be#1\ee{\begin{equation}#1\end{equation}}
\newcommand{\bea}{\begin{eqnarray}}
\newcommand{\eea}{\end{eqnarray}}
\newcommand{\beas}{\begin{eqnarray*}}
\newcommand{\eeas}{\end{eqnarray*}}

\newcommand{\tbA}{ \tilde{\bA}}
\newcommand{\tbB}{ \tilde{\bB}}
\newcommand{\tbC}{ \tilde{\bC}}

\newcommand{\cM}{{\mathcal{M}}}
\newcommand{\cS}{{\mathcal{S}}}
\newcommand{\cB}{{\mathcal{B}}}
\newcommand{\cX}{{\mathcal{X}}}
\newcommand{\cT}{{\mathcal{T}}}

\DeclareMathOperator*{\argmax}{arg\,max}

\title{Adaptive data-driven reduced-order models for reciprocal MIMO port-Hamiltonian dynamical systems}
\author{Anarzhan Abilgazy\thanks{Department of Mathematics, Southern Methodist University, Dallas, TX}, Mikhail Zaslavsky\thanks{Department of Mathematics, Southern Methodist University, Dallas, TX (\texttt{mzaslavskiy@smu.edu})}}
\date{}

\begin{document}

\maketitle

\begin{abstract}
Port-Hamiltonian (pH) systems provide a general framework for modeling a broad class of physical dynamical systems, including systems with energy storage, dissipation, and power exchange. Constructing accurate and computationally efficient data-driven reduced-order models (ROMs) for general pH systems is challenging, particularly when the underlying state-space realization is unavailable. In this work, we focus on an important subclass of pH systems satisfying reciprocity, a property arising naturally in a variety of wave, circuit, and network models. Rather than seeking to preserve the full pH structure, we exploit the additional symmetry provided by reciprocity to develop an efficient data-driven reduced-order modeling approach.

We show that every minimal realization of a reciprocal system is self-adjoint with respect to a generally indefinite inner product and can, up to a similarity transformation preserving the pH structure, be represented in a Krein space with a signature inner product. This observation motivates a Galerkin projection onto a tangential rational Krylov subspace with respect to the corresponding indefinite inner product. Although the resulting projection formally involves full-order state solutions and the unknown inner product, we show that its reduced mass and stiffness matrices can be constructed directly from transfer-function data and its frequency derivative through a Loewner-type formulation. The resulting reduced-order transfer function satisfies tangential Hermite interpolation conditions.

We then propose a greedy data-driven procedure that adaptively selects both interpolation frequencies and tangential directions. At every iteration, the frequency is chosen according to the largest transfer-function approximation error, while the tangential direction is selected from a dominant singular direction of the corresponding matrix-valued residual. Numerical experiments for wave propagation in lossy media demonstrate improved approximation efficiency compared with block adaptive interpolation. We further apply the method to spectral mass gauging of a three-dimensional propellant tank, where the data-driven ROM accurately recovers a large number of resonances and provides an accurate estimate of the tank volume through Weyl's law.    \end{abstract}
\maketitle

\section{Introduction}

Port-Hamiltonian (pH) systems constitute an important class of dynamical systems generalizing traditional Hamiltonian systems where the total energy is conserved to systems with energy dissipation. The pH formulation covers a broad range of physical models, including electrical circuits, mechanical and electromechanical systems, wave propagation in lossy media, and networked and distributed-parameter systems. In addition to providing a unified modeling framework, the pH structure encodes fundamental physical properties of the underlying system, e.g. its energy balance and, under standard assumptions, passivity and stability. For this reason, preserving the pH structure can be important when replacing a large-scale physical model by a reduced-order model (ROM).

Structure-preserving model reduction of pH systems has consequently received considerable attention; see, for example, \cite{GUGERCIN2012} and the references therein. The problem becomes more difficult, however, in a data-driven setting. Classical structure-preserving projection methods exploit the internal state-space representation and its structure, whereas in many applications only input-output measurements, such as samples of a transfer function, are available. Recovering a reduced model from such measurements while simultaneously retaining the structural properties of a general pH realization is nontrivial \cite{Benner2019IdentificationOP}.

A large and practically important subclass of pH systems possesses an additional property: reciprocity. In terms of inputs and outputs, reciprocity is expressed by symmetry of the MIMO transfer function,
\[
F(\omega)=F^T(\omega).
\]
This property occurs naturally in many physical systems. For example, wave propagation in lossy media governed by a symmetric differential operator with collocated sources and receivers is reciprocal, i.e. switching the source and receiver locations doesn't change the measured response. Similarly, reciprocal electrical networks provide another classical example. Both situations admit pH formulations and are considered explicitly in Section \ref{sec:ProblemSetup}. Thus, rather than identifying arbitrary pH systems from input-output data, in this work we focus on reciprocal pH systems and exploit the additional features brought by reciprocity. Note that our objective is not to preserve the entire port-Hamiltonian structure, but rather to exploit the additional symmetry provided by reciprocity to construct accurate and computationally efficient data-driven reduced-order models.

The key observation is that reciprocity, i.e. a symmetry of the measured transfer function, dictates a specific structure of dynamical systems itself. In particular, we show that any minimal realization of a reciprocal dynamical system is self-adjoint with respect to a generally indefinite inner product. Moreover, after an appropriate similarity transformation, this inner product can be represented by a signature matrix, while the pH structure of the realization is preserved. Reciprocal pH systems can therefore be treated as self-adjoint systems in a Krein space. This hidden self-adjointness is the structural property that enables the data-driven construction developed in this paper.

The Krein-space formulation naturally suggests a Galerkin projection \cite{Antoulas01asurvey} with respect to the corresponding indefinite inner product. Given frequency-domain solutions at a set of interpolation frequencies and a collection of tangential directions, the projection space is a tangential rational Krylov subspace \cite{GUGERCIN2012}. The resulting reduced mass and stiffness matrices are complex symmetric. More importantly, the symmetry of the projection yields a quadratic representation of the transfer-function approximation error. Consequently, interpolation of a tangential state direction implies Hermite tangential interpolation condition for the transfer function and its derivative rather than only tangential interpolation of the transfer function itself.

A drawback of the approach above for data-driven framework is that projection formally requires the full-order frequency-domain solutions as well as the indefinite inner product associated with reciprocity, neither of which is generally available from measurements. To resolve this issue, in this paper we developed a Loewner-type construction and we expressed the entries of the reduced mass and stiffness matrices entirely in terms of the transfer function and its frequency derivative evaluated at the interpolation frequencies. The projection-based ROM therefore admits an equivalent data-driven construction requiring no access to the full-order state-space realization or to the Krein inner product.

Once the reduced model can be constructed directly from data, an important remaining question is how the interpolation data should be selected. The accuracy of rational Krylov, and rational interpolation methods in general, depends strongly on the choice of interpolation frequencies. For frequency-domain problems with symmetric nonnegative-definite operators, a priori chosen optimized shifts based on Zolotarev approximation were proposed in \cite{knizhnerman2009optimal} and subsequently extended to evolutionary problems in \cite{druskin2009solution}. For Helmholtz problems in resonance domain, a priori sampling based on Weyl's law was proposed in \cite{Abilgazy}. These approaches to select interpolation points can be used within data-driven framework as well as in a model-driven setup where full-scale model and its solutions are accessible. For the latter only, adaptive rational Krylov methods were developed in \cite{druskin2010adaptive, druskin2011adaptive}, and adaptive tangential interpolation for MIMO systems was considered in \cite{druskin2014adaptive}. In turn, in the purely data-driven setting, the AAA algorithm \cite{aaa} provides an efficient adaptive rational approximation procedure, while algorithms for matrix-valued rational approximation, including block variants, were considered in \cite{blockaaa}. Vector fitting \cite{Gustavsen1999RationalAO} provides another widely used approach for constructing rational approximations from frequency-response data. However, all these data-driven approaches either treat every single element of transfer function independently or use full matrix-valued block version.

Here we develop an adaptive data-driven procedure specifically exploiting the tangential interpolation structure of reciprocal MIMO systems. Instead of interpolating the complete matrix-valued transfer function at every selected frequency, we adaptively select interpolation pair composed of the interpolation frequency and a tangential input-output direction. At each iteration, the new frequency is chosen where the current ROM has the largest error, and the tangential direction is chosen from a dominant singular direction of the matrix-valued residual at that frequency. The reflected frequency and conjugate tangential direction are added simultaneously to preserve the Schwarz reflection symmetry of the original frequency response. In this way, each new interpolation pair is directed toward a frequency region and input-output direction that are poorly represented by the current ROM. At the same time, the Krein-space projection provides Hermite interpolation of both the transfer function and its derivative in the selected tangential direction allowing to annihilate a large error of the previous ROM at selected interpolation pair.

The main contributions of this work are therefore as follows. First, we identify the state-space consequence of reciprocity that is central to our construction: a minimal reciprocal realization is self-adjoint in an indefinite inner product, and this representation is compatible with the pH structure. Second, exploiting this hidden self-adjointness, we construct a tangential rational Krylov ROM through a Galerkin projection in a Krein space and establish its tangential Hermite interpolation properties. Third, we show that the resulting reduced mass and stiffness matrices admit a fully data-driven Loewner-type representation in terms of transfer-function measurements and their frequency derivatives. Finally, we introduce an adaptive procedure that jointly selects interpolation frequencies and tangential directions and demonstrate numerically that, for the examples considered here, this provides a more efficient approximation in terms of ROM dimension than block adaptive interpolation.

We illustrate the approach with two numerical examples. The first considers wave propagation in a two-dimensional lossy medium and compares the proposed tangential adaptive method with block adaptive interpolation. Because block interpolation introduces a full set of input-output directions at every selected frequency, while the proposed method introduces only the dominant tangential direction, the latter can distribute a fixed ROM dimension over a larger number of frequencies. The numerical results demonstrate a substantial improvement in approximation accuracy for the considered examples. The second application concerns spectral mass gauging of a three-dimensional propellant tank. Here the data-driven ROM is used not only to approximate the measured frequency response but also to recover its resonances, which are subsequently employed to estimate the propellant volume through Weyl's law.

The remainder of the paper is organized as follows. Section \ref{sec:ProblemSetup} introduces reciprocal pH systems, presents representative physical examples, and establishes their Krein-space representation. Section \ref{sec:ProjROM} develops the projection-based model-driven ROM. Section \ref{sec:ddROM} derives its data-driven realization using the Loewner framework. Section \ref{sec:ROMProp} establishes interpolation and symmetry properties of the resulting model. Section \ref{sec:adaptiveROM} introduces the adaptive selection of interpolation frequencies and tangential directions. Numerical experiments are presented in Section \ref{sec:NumEx}, followed by concluding remarks in Section \ref{sec:conclusions}.

\section{Problem setup}
\label{sec:ProblemSetup}
	Consider multi-input-multi-output (MIMO) port-Hamiltonian dynamical system \cite{GUGERCIN2012} \bea\label{eq:pHDS1_td} 
	\frac{d\hx}{dt}=(\bJ-\bR)\bQ \hx+\bB \hu(t)\\\label{eq:pHDS2_td}
    \hy=\bB^T\bQ. \hx\eea
Here $\hx(t)\in\RR^N$ is state vector, $\hu(t)\in\RR^m$ is input, $\hy(t)\in\RR^m$ is output, $\bB\in\RR^{N\times m}$ and $\bC=\bB^T\bQ$ are input and output matrices, $\bJ=-\bJ^T\in\RR^{N\times N}$ is skew-symmetric structure matrix and $0\le\bR=\bR^T\in\RR^{N\times N}$ is symmetric non-negative dissipation matrix,
$\bQ=\bQ^T>0$ is symmetric energy matrix that defines Hamiltonian $H(\hx)=\frac{1}{2}\hx^T\bQ \hx$. Below we shall relate the dynamical system (\ref{eq:pHDS1_td})-(\ref{eq:pHDS2_td}) to a triplet $[\bA,\bB,\bC]$ where $\bA=(\bJ-\bR)\bQ$ is system matrix. Also, here and below $^T$ means real transpose, even when applied to complex matrices.

    The Fourier frequency-domain counterpart of the system (\ref{eq:pHDS1_td}), (\ref{eq:pHDS2_td}) has a form
    \bea\label{eq:pHDS1} 
	i\omega x-(\bJ-\bR)\bQ x=\bB u(\omega)\\\label{eq:pHDS2}
    y=\bB^T\bQ x.\eea Here $x$, $y$ and $u$ are Fourier images of $\hx$, $\hy$ and $\hu$, respectively. 
We define frequency-domain MIMO transfer function as 
    \be F(\omega)=\bB^T\bQ(i\omega \bI-(\bJ-\bR)\bQ)^{-1}\bB\in\CC^{m\times m}\ee where $\bI\in\RR^{N\times N}$ is the identity matrix. Note that $i$-th column of $F(\omega)$ is output $y$ of dynamical system (\ref{eq:pHDS1}), (\ref{eq:pHDS2}) corresponding to input being $i$-th unit vector: $u(\omega)=e_i$.

Clearly, the block of solutions $\bx=(i\omega \bI-(\bJ-\bR)\bQ)^{-1}\bB\in\CC^{N\times m}$ as well the transfer function $F$ satisfy Schwarz reflection symmetry condition
\be \bx(-\omega)=\bar{\bx}(\omega),~F(-\omega)=\bar{F}(\omega)\ee where $\bar{\cdot}$ denote complex conjugation. This property guarantees that the time-domain counterpart $\hat{\bx}(t)$ of $\bx(\omega)$ is real-valued.

In this paper we focus on a reciprocal subclass of port-Hamiltonian systems defined as follows:
\begin{definition}
MIMO dynamical system is called reciprocal if its transfer function is symmetric $F(\omega)=F^T(\omega)$.
\end{definition}

To illustrate the concept, below we provide two examples of reciprocal port-Hamiltonian systems
\begin{examp}\label{ex1}
Consider damped Helmholtz equation in the second-order form \be\label{dampedhelm} -\Delta u^j(\vec{r})+i\omega R u^j(\vec{r})-\omega^2 u^j(\vec{r})=i\omega f^j(\vec{r}),~\vec{r}\in\RR^d,~d\ge 1\ee with symmetric non-negative damping operator $R\ge 0$ as well as multiple inputs $f^j(\vec{r}),~j=1,\ldots, m$ and outputs given by $y^{ij}(\omega)=\int{u^j(\vec{r})f^i(\vec{r})d\vec{r}}$. Its finite-difference discretization has a form

\be\label{2ndord} L^T_hL_hu^j_h+i\omega R_hu^j_h-\omega^2u^j_h=i\omega f^j_h. \ee
with outputs $y^{ij}_h=(f^i)^T_hu^j_h$. Here $u^j_h\in\CC^{N_1}$, $f^j_h\in\RR^{N_1}$, $\RR^{N_1\times N_1}\ni R_h=R^T_h\ge 0$, $L_h\in\RR^{N_2\times N_1}$ is discretization of the gradient operator $\nabla$, and $L^T_h$ is its adjoint discretizing the minus divergence operator $-\nabla\cdot$. We note that all the operators in (\ref{2ndord}) are symmetric and, therefore, the system is reciprocal $y^{ij}=y^{ji}$. Then, introducing an auxiliary variable $v^j_h=\frac{1}{i\omega}L_hu^j_h\in\CC^{N_2}$, from (\ref{2ndord}) we obtain first-order formulation

    \begin{equation}\label{1stord}
    i\omega\begin{pmatrix}u^j_h\\v^j_h \end{pmatrix}-\left(\begin{pmatrix}0&-L^T_h\\L_h&0 \end{pmatrix}-\begin{pmatrix}R_h&0\\0&0 \end{pmatrix}\right)\begin{pmatrix}u^j_h\\v^j_h \end{pmatrix}=\begin{pmatrix}f^j_h\\0 \end{pmatrix}
\end{equation}
with outputs given by $$y^{ij}=\begin{pmatrix}f^i_h\\0 \end{pmatrix}^T\begin{pmatrix}u^j_h\\v^j_h \end{pmatrix}.$$ Collecting all the inputs $f^j_h,~j=1,\ldots, m$ in one matrix $$\bB=\left(\begin{pmatrix}f^1_h\\0 \end{pmatrix},\ldots,\begin{pmatrix}f^m_h\\0 \end{pmatrix}\right)\in\RR^{N\times m}$$ with $N=N_1+N_2$ and all the solutions for all inputs in another one $$\bx=\left(\begin{pmatrix}u^1_h\\v^1_h \end{pmatrix},\ldots,\begin{pmatrix}u^m_h\\v^m_h \end{pmatrix}\right)\in\RR^{N\times m},$$ we can rewrite (\ref{1stord}) as 
    \begin{equation}\label{1stord_block}
    i\omega\bx-\left(\begin{pmatrix}0&-L^T_h\\L_h&0 \end{pmatrix}-\begin{pmatrix}R_h&0\\0&0 \end{pmatrix}\right)\bx=\bB
\end{equation}
where outputs are given by $y=\bB^T\bx\in\CC^{m\times m}$. Clearly, this reciprocal system has MIMO frequency-domain port-Hamiltonian form (\ref{eq:pHDS1})-(\ref{eq:pHDS2}) with structure, dissipation and energy matrices given by
    \begin{equation}
    \bJ=\begin{pmatrix}0&-L^T_h\\L_h&0 \end{pmatrix},~\bR=\begin{pmatrix}R_h&0\\0&0 \end{pmatrix},~\bQ=\bI,
\end{equation} respectively.
\end{examp}
\begin{examp}\label{ex2}
Consider a system of ODEs for $m$ electric charge evolutions in RCL electric circuit due to $m$ independent inputs. Matrix $q\in \RR^{N_1\times m}$ storing those evolutions columnwise satisfies 
\be\label{rcl2ndord} Lq''+Rq'+Cq=\bB_1\ee where columns of $\bB_1\in\RR^{N_1\times m}$ correspond to $m$ inputs. Here $0<L=L^T\in \RR^{N_1\times N_1}$, $0\le R=R^T\in \RR^{N_1\times N_1}$ and $0<C=C^T\in \RR^{N_1\times N_1}$ are positive symmetric inductance, non-negative symmetric resistance  and positive symmetric capacitance  matrices, respectively. Denote reciprocal outputs $y=\bB_1^Tq'$. Introducing an auxiliary variable $p=Lq'$ from (\ref{rcl2ndord}) we obtain 
    \begin{equation}\label{rcl1stord}\begin{pmatrix}p\\q \end{pmatrix}'=
    \left(\begin{pmatrix}0&-\bI\\ \bI&0 \end{pmatrix}-\begin{pmatrix} R &0\\0&0 \end{pmatrix}\right)\begin{pmatrix}L^{-1}&0\\0&C \end{pmatrix}\begin{pmatrix}p\\q \end{pmatrix}+\begin{pmatrix}\bB_1\\0 \end{pmatrix}
\end{equation}
with the output given by $y=\bB_1^Tq'=\begin{pmatrix}\bB_1 \\ 0 \end{pmatrix}^T\begin{pmatrix}L^{-1}&0\\0& C \end{pmatrix}\begin{pmatrix}p\\q \end{pmatrix}$. Clearly, this reciprocal system has MIMO port-Hamiltonian form (\ref{eq:pHDS1_td})-(\ref{eq:pHDS2_td}) with structure, dissipation, energy and input matrices     \begin{equation}
    \bJ=\begin{pmatrix}0&-\bI\\ \bI&0 \end{pmatrix},~\bR=\begin{pmatrix}R&0\\0&0 \end{pmatrix},~\bQ=\begin{pmatrix}L^{-1}&0\\0& C\end{pmatrix},~\bB=\begin{pmatrix}\bB_1\\0 \end{pmatrix},
\end{equation} respectively.

\end{examp}
We note that system matrices in Examples \ref{ex1} and \ref{ex2} are so-called Krein-self-adjoint, i.e. they are symmetric with respect to indefinite inner products. Indeed, these inner products are induced by indefinite matrices $\Sigma=\begin{pmatrix}\bI_{N_1}&0\\0&-\bI_{N_2} \end{pmatrix}$ and $\Sigma=\begin{pmatrix}L^{-1}&0\\0&-C \end{pmatrix}$ for Example \ref{ex1} and \ref{ex2}, respectively. Moreover, output matrices $\bB^T\bQ$ in both examples are equal to $\bB^T\Sigma$. The following proposition shows that it is inherent property of all minimal reciprocal systems. Moreover, up to similarity transform preserving the port-Hamiltonian structure, the reciprocal systems are symmetric with respect to the inner products induced by signature matrices, i.e. diagonal matrices with either $1$'s or $-1$'s on diagonal. 

\begin{prop}
\label{prop_sym}
    Let $[\bA,\bB,\bC]$ be minimal realization  of dynamical system with reciprocal transfer function $F(\omega)=\bC(i\omega-\bA)^{-1}\bB=F^T(\omega)$.
    \begin{enumerate}
        \item Then there exists a symmetric matrix $\Sigma=\Sigma^T$ such that $\bA$ is $\Sigma$-symmetric: $\Sigma \bA=\bA^T\Sigma$ and $\bC=\bB^T\Sigma$
        \item There exists an equivalent system $[\tbA,\tbB,\tbC]$ that is symmetric with respect to a signature matrix $\Sigma=\begin{pmatrix}\bI_{N_1}&0\\0&-\bI_{N_2} \end{pmatrix}$ for some values $N_1$ and $N_2$
        \item If $[\bA,\bB,\bC]$ is port-Hamiltonian then $[\tbA,\tbB,\tbC]$ is port-Hamiltonian as well
    \end{enumerate}
\end{prop}
\begin{proof}\nonumber
    \begin{enumerate}
    \item We have $\bB^T(i\omega-\bA^T)^{-1}\bC^T=F^T(\omega)=F(\omega)$. Therefore, $[\bA^T,\bC^T,\bB^T]$ is another realization of the dynamical system producing the same transfer function $F(\omega)$. Since $[\bA,\bB,\bC]$ is minimal, $[\bA^T,\bC^T,\bB^T]$ differs from $[\bA,\bB,\bC]$ by similarity transform, i.e. there exists such $\Sigma$ that $\bC=\bB^T\Sigma$, $\bB=\Sigma^{-1}\bC^T $ and $\bA^T=\Sigma\bA\Sigma^{-1} $. To show that $\Sigma$ is symmetric, note that $\bA=\Sigma^{-T}\bA^T\Sigma^T=\Sigma^{-T}\Sigma\bA\Sigma^{-1}\Sigma^T$, or, equivalently, $\bA\Sigma^{-T}\Sigma=\Sigma^{-T}\Sigma\bA$. Next, from $\bC^T=\Sigma\bB$ we have $\Sigma^{-T}\Sigma\bB=\Sigma^{-T}\bC^T$. But from $\bC=\bB^T\Sigma$ we have $\Sigma^{-T}\bC^T=\bB$. Therefore, $\Sigma^{-T}\Sigma\bB=\bB$. Then, since $\Sigma^{-T}\Sigma$ commutes with $\bA$, for $k\ge 0$ we have: $\Sigma^{-T}\Sigma\bA^k\bB=\bA^k\Sigma^{-T}\Sigma\bB=\bA^k\bB$. Finally, since $[\bA,\bB]$ is controllable, vectors $\{\bA^k\bB\}$ span the entire space, so $\Sigma^{-T}\Sigma=\bI$, i.e. $\Sigma=\Sigma^T$
    \item By Sylvester law of inertia, there exists $\bS$ such that $$\Sigma=\bS^T\begin{pmatrix}\bI_{N_1}&0\\0&-\bI_{N_2} \end{pmatrix}\bS$$ where $N_1$ and $N_2$ are positive and negative indices of inertia, respectively. Then we can set $\tbA=\bS\bA\bS^{-1}$, $\tbB=\bS\bB$ and $\tbC=\bC\bS^{-1}=\tbB^T\begin{pmatrix}\bI_{N_1}&0\\0&-\bI_{N_2} \end{pmatrix}$
    \item Indeed, if $\bA=(\bJ-\bR)\bQ$ then $\tbA=\bS(\bJ-\bR)\bQ\bS^{-1}=(\bS\bJ\bS^T-\bS\bR\bS^T)\bS^{-T}\bQ\bS^{-1}$ that corresponds to structure, dissipation and energy matrices given by $\bS\bJ\bS^T$, $\bS\bR\bS^T$ and $\bS^{-T}\bQ\bS^{-1}$, respectively. Finally, for input matrix $\bS\bB$ the output matrix is $$\bB^T\bQ\bS^{-1}=(\bB^T\bS^T)(\bS^{-T}\bQ\bS^{-1})$$ which is consistent with port-Hamiltonian property.
    \end{enumerate}
\end{proof}
    This proposition allows to handle reciprocal dynamical systems in selfadjoint way in some Krein space equipped with indefinite inner product. In the next section we shall exploit that fact to construct reduced-order models (ROMs) of reciprocal port-Hamiltonian systems. We note that finding such inner product requires operations with full system $[\bA,\bB,\bC]$ that may be cost prohibitive for large-scale scenarios. To avoid that, we shall approach this problem within data-driven framework and construct ROM using transfer function only via Loewner framework, as will be discussed in Section \ref{sec:ddROM}.

\section{Projection-Based Model-Driven ROM}
\label{sec:ProjROM}
Let us consider the equation
\begin{eqnarray}
    i\omega \bx -\bA\bx  = \bB 
\end{eqnarray}
with $\bA=(\bJ - \bR)\bQ$.
Let $\bx_q \in\CC^{N\times m} $ be its solutions corresponding to  $\omega=\omega_q,~q=1,...,n$ and denote $F_q=F(\omega_q)$. We also pick some vectors $r_q\in\CC^m, q=1,\ldots n$ that are typically called (right) tangential vectors. Proposition \ref{prop_int} explains the reason for this name. We shall assume that there are no coinciding pairs $(\omega_q,r_q)$. A model-driven reduced-order model (ROM) can then be constructed using a Galerkin formulation by projecting the problem onto the rational Krylov subspace $\mathcal{V} = colspan\{ \bx_1 r_1,..., \bx_n r_n\}$ with respect to indefinite inner product induced by matrix $\Sigma$ from Proposition \ref{prop_sym}. Denoting $\bV=[\bx_1 r_1,...,\bx_n r_n]$  $\in\CC^{N\times n}$, we obtain the following.
 \begin{eqnarray} \label{eq:projsyst}
     (i\omega \cM-\cS) c = \cB.
\end{eqnarray} 
The matrix $\cM=\cM^T=\bV^T\Sigma\bV$, a complex-symmetric indefinite mass matrix with entries
\begin{equation}\label{massmtr}\cM_{qp}=r^T_q \bx^T_q  \Sigma \bx_p r_p;  ~q,p=1,\ldots , n\end{equation}
Similarly, the second term defines the matrix $\cS=\cS^T=\bV^T\Sigma\bA\bV$, a complex-symmetric indefinite stiffness matrix with entries
 \begin{equation} 
\label{stifmtr}
\cS_{qp}=r^T_q\bx^T_q \Sigma \bA \bx_p r_p;  ~q,p=1,\ldots , n.\end{equation}
and entries of $\cB=\bV^T\Sigma\bB \in \mathbb{C}^{n \times m}$ are given by
\begin{eqnarray*}
\cB_q = r^T_q\bx^T_q\Sigma\bB.
\end{eqnarray*}
The ROM approximation of the solution $\bx$ is given by
\begin{eqnarray}
\label{eq:galsol_2d}
    \bx(\omega) \approx \cX(\omega) = \bV c = \bV (i\omega \cM-\cS)^{-1}  \cB.
\end{eqnarray}

\noindent Our reduced order MIMO transfer function is given by
\begin{eqnarray}\label{F_rom}
F^{ROM}(\omega) = \bB^T\Sigma\cX=\cB^T c=\cB^T (i\omega \cM-\cS)^{-1}  \cB.
\end{eqnarray}
It is reciprocal because $\cM$ and $\cS$ are complex-symmetric, i.e.
\be\nonumber (F^{ROM}(\omega))^T = (\cB^T (i\omega \cM-\cS)^{-1}  \cB)^T=\cB^T (i\omega \cM^T-\cS^T)^{-1}  \cB=\cB^T (i\omega \cM-\cS)^{-1}  \cB=F^{ROM}(\omega).\ee
Assuming that matrix pencil $(\cS;\cM)$ has a complete $\cM$-orthogonal basis composed of eigenvectors $z_i$ $$\cS z_i=\theta_i\cM z_i, i=1,\ldots ,n$$
$$z^T_i\cM z_j=\delta_{ij}$$
we obtain that $F^{ROM}$ can be represented in terms of symmetric matrix rational function
\be\label{ratfunc} F^{ROM}(\omega)=(F^{ROM}(\omega))^T=\sum^n_{i=1}\frac{\cB^Tz_iz^T_i\cB}{i\omega-\theta_i}.\ee

\section{Data-driven framework}
\label{sec:ddROM}
In multiple applications access to ROM stiffness and mass matrices using formulas (\ref{stifmtr}) and (\ref{massmtr}) is either restricted or computationally expensive. For example, in inverse scattering and remote sensing the operator $\bA$ is unknown, and the goal is to find it in some parametric form from the measured data given by transfer function. That said, inner product $\Sigma$ as well as solutions $\bx_p$ required to compute stiffness and mass matrices in ROM (\ref{stifmtr}) and (\ref{massmtr}) are not accessible either. Also, when $\bA$ is dense and size of ROM $n$ is large, the computational cost of computing multiple matrix-vector products in (\ref{stifmtr}) and (\ref{massmtr}) can be prohibitively high. However, the mass and stiffness matrices can still be computed directly using transfer function only via Loewner framework. Indeed, for $\omega = \omega_p$ the projected system can be written as

\begin{equation}
\label{eq:krylov}
    i\omega_p r^T_q x^T_q \Sigma  x_p r_p-r^T_qx^T_q \Sigma \bA x_p r_p =r^T_q x^T_q \Sigma \bB r_p
\end{equation}
Due to symmetry of transfer function, the right-hand side of \eqref{eq:krylov} is given by 
\begin{equation}\label{D_rom}
r^T_q x^T_q  \Sigma\bB r_p = r^T_q F^T_q r_p= r^T_q F_q r_p = r^T_q F_q( \omega) r_p\in\CC~q=1,\ldots , n.\end{equation}
Therefore, we obtain 
\begin{equation}
\label{omegap} 
i\omega_p \cM_{qp}-\cS_{qp}= r^T_q F_q r_p
\end{equation}
or, switching indexes, 
\begin{equation}
\label{omegaT}
    i\omega_p \cM_{pq}-\cS_{pq}= r^T_p F_q  r_q.
\end{equation}

\noindent Similarly, for $\omega = \omega_q$ we have 
\begin{equation}
\label{omegaq}
i\omega_q \cM_{pq}-\cS_{pq} = r^{T}_{p} F_p r_q.
\end{equation}
Subtracting \eqref{omegaT} from  \eqref{omegaq}, we obtain
\be
    \cM_{pq} (i\omega_q - i\omega_p) = r^T_p F_p r_q - r^T_p F_qr_q
\ee
Hence, elements of the mass matrix can be expressed in a purely data-driven form as
\begin{equation}\label{massmtrdd} 
\cM_{pq}= r^T_px^T_p  \Sigma  x_qr_q= r^T_p \frac{F_p-F_q}{i\omega_q-i\omega_p} r_q. \end{equation}
Similarly, expressing elements $\cS_{pq}$ of the stiffness matrix from (\ref{omegaq}) and (\ref{omegaT}), we obtain 
\begin{eqnarray} \label{S}
       \cS_{pq} = -r^T_p\frac{\omega_p F_p -\omega_q F_q}{\omega_p-\omega_q}r_q. 
\end{eqnarray}
The case $p = q$ follows from (\ref{massmtrdd}) and (\ref{S}) by formally passing $\omega_q\rightarrow\omega_p$:
\begin{eqnarray}\label{diagM}
     \cM_{pp} & = & i r^T_p  \left.\frac{dF}{d\omega}\right|_{\omega=\omega_p} r_p,\\  
     \cS_{pp} & = & -r^T_p (F_p  + \omega_p \left.\frac{dF}{d\omega}\right|_{\omega=\omega_p}) r_p. \label{diagS}
\end{eqnarray}

\section{Properties of the obtained ROM}
\label{sec:ROMProp}
The following proposition establishes the interpolation properties of the constructed ROM.
\begin{prop} \label{interpprop} Let matrix $i\omega_p\cM-\cS$ be invertible. Then the following interpolation conditions hold
 \begin{enumerate}
 \item $\bx_pr_p=\cX(\omega_p)r_p$
 \item $F^{ROM}(\omega_p)r_p=F_pr_p$, $r^T_pF^{ROM}(\omega_p)=r^T_pF_p$
 \item $r^T_p\left.\frac{dF^{ROM}}{d\omega}\right|_{\omega=\omega_p}r_p=r^T_p\left.\frac{dF}{d\omega}\right|_{\omega=\omega_p}r_p$
 \end{enumerate}
    
\label{prop_int}
\end{prop}
\begin{proof}
\begin{enumerate}
\item Since for $p=1,\ldots n$ we have   $i\omega_p \bx_p -\bA \bx_p  = \bB$, multiplying it by $r_p$ from the right and by $r_q^T\bx^T_q\Sigma$ for arbitrary $q$ from the left, we obtain $i\omega_p\cM_{qp}-\cS_{qp}=\cB_qr_p$. Therefore, solution $c$ of  (\ref{eq:projsyst}) for $\omega=\omega_p$ satisfies $cr_p=e_p$ where $e_p$ is $p$-th unit vector. Due to assumptions of this Proposition, this solution is unique. Hence, since $\cX=\bV c$, we have \be\label{solint} \cX(\omega_p)r_p=\bV cr_p=\bV e_p=\bx_pr_p\ee
\item First note that from (\ref{eq:projsyst}) $c$ satisfies $\bV^T\Sigma\bB-(i\omega \bV^T\Sigma\bV-\bV^T\Sigma\bA\bV)c=0$, or, equivalently, $\bV^T\Sigma (\bB-(i\omega I-\bA)\cX)=0$. Therefore, taking into account that $\cB^Tc=\bB^T\Sigma\cX$, $\bB=(i\omega\bI-\bA)\bx$ as well as $\Sigma A=A^T\Sigma$, we have 
\bea\label{tfquadr}\nonumber F(\omega)-F^{ROM}(\omega)=\bB^T\Sigma\bx-\cB^Tc=\\ \nonumber \bB^T\Sigma\bx-\bB^T\Sigma\cX-c^T\bV^T\Sigma (\bB-(i\omega \bI-\bA)\cX)=\\ \nonumber  \bB^T\Sigma(\bx-\cX)-\cX^T\Sigma(i\omega\bI-\bA)(\bx-\cX)=\\ \nonumber  (\bx^T(i\omega\bI-\bA)^T\Sigma-\cX^T\Sigma(i\omega\bI-\bA))(\bx-\cX)=\\ (\bx-\cX)^T\Sigma(i\omega\bI-\bA)(\bx-\cX).\eea Hence, multiplying (\ref{tfquadr}) by $r_p$ from the right and evaluating it for $\omega=\omega_p$, we obtain \be \label{tfquadrrp} F(\omega_p)r_p-F^{ROM}(\omega_p)r_p=(\bx-\cX)^T\Sigma(i\omega\bI-\bA)(\bx r_p-\cX r_p)=0\ee because of (\ref{solint}). Then, taking the reciprocity of $F(\omega)$ and $F^{ROM}(\omega)$ into account, we have \bea\nonumber r^T_pF(\omega_p)-r^T_pF^{ROM}(\omega_p)=(F^T(\omega_p)r_p-(F^{ROM})^T(\omega_p)r_p)^T=\\ (F(\omega_p)r_p-F^{ROM}(\omega_p)r_p)^T=0.\eea
\item Let's multiply (\ref{tfquadr}) by $r_p$ from the left and from the right and differentiate the result with respect to $\omega$. Using (\ref{solint}) finalizes the proof.
\end{enumerate}
\end{proof}
If assumptions of Proposition \ref{prop_int} hold then the transfer function is said to satisfy tangential interpolation conditions at $\omega_q$ with (right) tangential vectors $r_q$ (hence, the name for vectors $r_q$). The pair $(\omega_q,r_q)$ is called interpolation pair, and a set of all such pairs $I=\{(\omega_q,r_q)\}^n_{q=1}$ is called interpolation set.

The choice of interpolation pairs $(\omega_q,r_q),~q=1,\ldots ,n$ is crucial for constructing ROM that is accurate and preserves the properties of the original full-scale problem. To satisfy Schwarz reflection symmetry condition, every interpolation pair $(\omega_q,r_q)$ should enter the interpolation set together with $(-\omega_q,\bar{r}_q)$. In particular, the following proposition holds 

\begin{prop} \label{schwarzprop}
If interpolation set $I$ satisfies the property that for any pair $(\omega_q,r_q)\in I$ a pair $(-\omega_q,\bar{r}_q)\in I$ then the obtained ROM satisfies Schwarz reflection symmetry condition\be \cX(-\omega)=\bar{\cX}(\omega),~F^{ROM}(-\omega)=\bar{F}^{ROM}(\omega)\ee
\end{prop}
\begin{proof}
Note that under assumptions of this proposition the columns of $\bar{V}$ span the same subspace as columns of $V$. Therefore, there exists such invertible matrix $\cT\in\CC^{n\times n}$ that $\bar{V}=V\cT$. From (\ref{eq:galsol_2d}) we have \bea\nonumber\bar{\cX}(\omega)=\bar{\bV}\bar{c}=\bar{\bV}(-i\omega\bar{\cM}-\bar{\cS})^{-1}\bar{\cB}=\\ \nonumber \bar{\bV}(-i\omega\bar{\bV}^T\Sigma\bar{\bV}-\bar{\bV}^T\Sigma\bA\bar{\bV})^{-1}\bar{\bV}^T\Sigma\bB=\\ \nonumber \bV \cT(-i\omega \cT^T\bV^T\Sigma\bV \cT-\cT^T\bV^T\Sigma\bA\bV \cT)^{-1}\cT^T\bV^T\Sigma\bB=\\ \nonumber\bV (-i\omega\bV^T\Sigma\bV-\bV^T\Sigma\bA\bV )^{-1}\bV^T\Sigma\bB=\\ \bV (-i\omega\cM -\cS)^{-1}\cB=\cX(-\omega).\eea Then $$\bar{F}^{ROM}(\omega)=\bB^T\bar{\cX}(\omega)=\bB^T\cX(-\omega)=F^{ROM}(-\omega).$$
\end{proof}
\begin{remark}
    In most of the cases, the obtained projection-based reduced-order model is not structure-preserving, i.e. it doesn't inherit port-Hamiltonian structure of the original problem. As such, it may lose passivity and even stability property of the full-order system. One of the exceptions here corresponds to Hamiltonian case of Example \ref{ex1}, i.e. where loss term $\bR$ is vanishing. One of the approaches to preserve port-Hamiltonian structure was developed in \cite{Benner2019IdentificationOP} and is based on interpolation at specific points, i.e. at spectral zeros. While it allows to inherit passivity of the original system, the accuracy of the obtained ROM can be far from optimal. In turn, in our paper we focus selecting interpolation points adaptively by minimizing data misfit, as will be discussed below. 
\end{remark}

\section{Adaptive approach to select interpolation pairs}
\label{sec:adaptiveROM}
Let's consider the problem of optimizing the choice of interpolation pairs $(\omega_q,r_q)$ for accuracy of transfer function $F(\omega)$ approximation by $F^{ROM}(\omega)$ on frequency interval $[\omega_{min};\omega_{max}]$, For the single-input-single-output (SISO) scenario and symmetric positive-definite operator $\bA$ model-driven framework for this problem was discussed in \cite{knizhnerman2009optimal} where authors developed an algorithm to pick interpolation frequencies based on Zolotarev approximant. The approach was further extended in \cite{druskin2009solution} to handle time-domain problems. Such choice has been shown to be optimal for the case of uniform spectral measure of operator $\bA$. To improve the accuracy in the case of highly non-uniform spectral measure, authors of \cite{druskin2010adaptive} developed model-driven adaptive approach based on greedy algorithm. It was further generalized to handle MIMO scenarios with symmetric positive-definite $\bA$ in \cite{druskin2014adaptive} and SISO scenarios with non-symmetric $\bA$ in \cite{druskin2011adaptive}. Data-driven selection of interpolation frequencies can be performed via Adaptive Antoulas-Andersen (AAA) approach. SISO scenario was considered in \cite{aaa}, and block-AAA for MIMO problems was later proposed in \cite{blockaaa}. Another approach to construct rational approximant was proposed in \cite{Gustavsen1999RationalAO}. It allows to construct passive rational function in an iterative way, however it is not tied to selection of interpolation points, and its convergence heavily depends on efficient choice of initial guess. In this paper we propose an adaptive approach to choose not just interpolation points $\omega_q,~q=1,\ldots ,n$ but also tangential vectors $r_q,~q=1,\ldots n$ in order to construct rational approximant satisfying tangential interpolation conditions. In particular, we choose interpolation pair $(\omega_q,r_q)$ at each step in a greedy way, i.e. first pick $\omega_q$ maximizing the error $\|F(\omega)-F^{ROM}(\omega)\|$ for $\omega\in[\omega_{min};\omega_{max}]$ and then pick $r_q$ maximizing $\|(F(\omega_q)-F^{ROM}(\omega_q))r_q\|$. As a result, the obtained interpolation pair maximizes $$\max_{\omega,\|r\|_2=1}{\|(F(\omega)-F^{ROM}(\omega))r\|}$$ Thanks to Hermite tangential interpolation properties established in Proposition \ref{prop_int}, at the next step  we shall have $(F(\omega_q)-F^{ROM}(\omega_q))r_q=0$ and $r^T_q\left.\frac{d(F-F^{ROM})}{d\omega}\right|_{\omega=\omega_q}r_q=0$. In order to preserve Schwarz reflection symmetry, we add $(-\omega_q;\bar{r}_q)$ as well. In our numerical examples we shall show that it allows to obtain a superior accuracy compared to interpolation of the full matrix transfer function via block algorithm.  
\begin{algorithm}[H]
\caption{Adaptive Rational Tangential Data Fitting }
\label{alg:adpol}
\begin{algorithmic}[1]
\State Identify the data block $F_m$ with the maximum 2-norm and its corresponding frequency $\omega_m$. Set 
$\omega_1 = \omega_m$ and $\omega_2 = -\omega_1$.
\State Perform the SVD of $F_m$ and set the first tangential vector $r_1$ equal to the leading right singular vector $v_1$. Then set $r_2 = \Bar{r}_1$
\For{$n = 2, 4, 6, \ldots$ until convergence}
    \State Compute the matrices $\cS$ and $\cM$ via \eqref{massmtrdd}, \eqref{S}, \eqref{diagM}, \eqref{diagS} for interpolation frequencies $\omega_1, \ldots, \omega_n$.
    
    \State Compute the ROM data $F^{\mathrm{ROM}}$ for $\omega \in [\omega_{\min}, \omega_{\max}]$ via \eqref{F_rom}.
    \State Compute the next interpolation frequency by maximizing the approximation error:
    $\displaystyle
    \omega_{n+1} =
    \argmax_{\omega \in [\omega_{\min}, \omega_{\max}]}
    \left\|F^{\mathrm{ROM}}(\omega)-F(\omega)\right\|_2.
    $
    If multiple frequencies attain the maximum, select any one of them.
    \State Set $r_{n+1}$ to the leading right singular vector of $[F^{\mathrm{ROM}}(\omega_{n+1}) - F(\omega_{n+1})]$.
    \State Set $\omega_{n+2} = -\omega_{n+1}$ and $r_{n+2} = \Bar{r}_{n+1}$.
  
\EndFor
\end{algorithmic}
\end{algorithm}
\noindent

\begin{remark}
 After step 4, the eigendecomposition of mass matrix $\cM$ can be used to remove numerically insignificant and noisy modes and noise from the ROM. In particular, all the eigenvalues of $\cM$ with absolute values below the prescribed relative tolerance $\varepsilon$ are discarded. The ROM matrices and data are then projected onto this eigenbasis to improve numerical stability in the subsequent rational approximation.

\end{remark}
\noindent
Note that block-AAA-type approaches pick interpolation points $\omega_q,~q=1,\ldots$ in a similar greedy way, however then interpolation conditions for the full transfer function $F(\omega)\in\CC^{m\times m}$ are enforced by adding the entire block $\bI_m$ of $m$ tangential vectors simultaneously. Also, AAA employs baricentric representation for matrix rational interpolant and obtains its remaining parts via minimizing $L_2$ data misfit rather than interpolating $\frac{dF}{d\omega}$. In turn, our ROM enforces Hermite tangential interpolation conditions outlined in Proposition \ref{interpprop}. That allows us to focus on the parts of the transfer function that are hard to approximate and, consequently, to produce more accurate approximant, though it is applicable to interpolation of transfer functions of dynamical system (\ref{eq:pHDS1_td},\ref{eq:pHDS2_td}) only, i.e. when Loewner framework can be applied.

\section{Numerical examples}
\label{sec:NumEx}
\subsection{2D example in lossy media}
\begin{figure}[t]
\centering
\begin{tabular}{c @{\hspace{-0.4cm}}c}
Loss & $\tfrac{1}{2} \nabla  \log(\sigma)$ \\
\includegraphics[height=4.65cm]{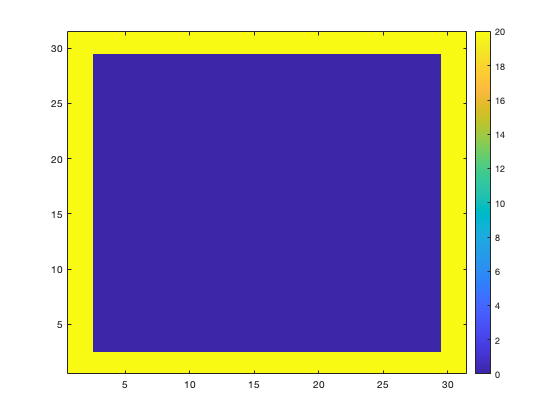} &
\includegraphics[height=4.65cm]{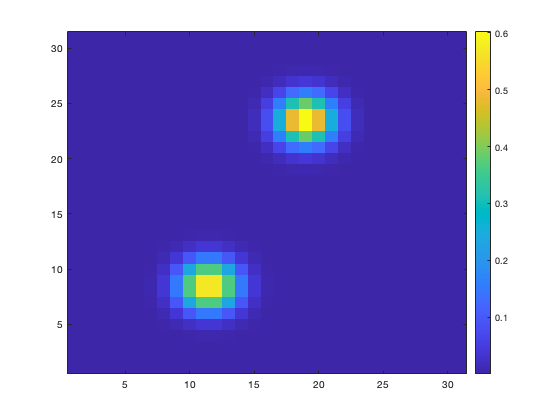}
\end{tabular}
\caption{{\bf Left:} The loss function is defined as a piecewise-constant function on the unit square. It is set to zero in the interior region \(0.05\leq x,y\leq0.95\), while taking the value \(20\) in the surrounding boundary layer. \\
{\bf Right:} The function involving the impedance  ($\tfrac{1}{2} \nabla \log(\sigma)$, where $\sigma$ is impedance) is modeled as the superposition of two localized Gaussian profiles. The Gaussian components are centered at \((0.6,0.35)\) and \((0.75,0.25)\), respectively, with standard deviation \(0.05\) in each spatial direction. The resulting coefficient is scaled by a factor of 0.01.
}
\label{fig: loss}
\end{figure}

We consider port-Hamiltonian system (\ref{1stord_block}) obtained from the finite-difference discretization of 2D damped Helmholtz equation on uniform Cartesian grid with $n_x=30$ and $n_y=30$ grid cells along $x$ and $y$, respectively. That results in dynamical system sizes 
$N_1=(n_x+1)(n_y+1)=961, ~N_2=n_x(n_y+1)+(n_x+1)n_y=1860,~N=N_1+N_2=2821.$ The medium parameters are shown  in Figure \ref{fig: loss}. We consider frequencies in the range
$\omega \in  [0.1,30]$, with \(N_{\omega}=500\) uniformly spaced frequency points in each interval. We compared $m=8$ and $m=16$ sources (distributed 2 per side and 4 per side, respectively).

\begin{figure}[t]
\centering
\begin{tabular}{c @{\hspace{-0.4cm}} c}
2 sources per side & 4 sources per side \\
\includegraphics[height=4.6cm]{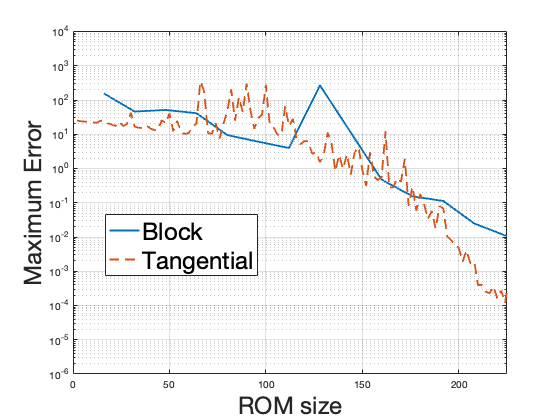} &
\includegraphics[height=4.6cm]{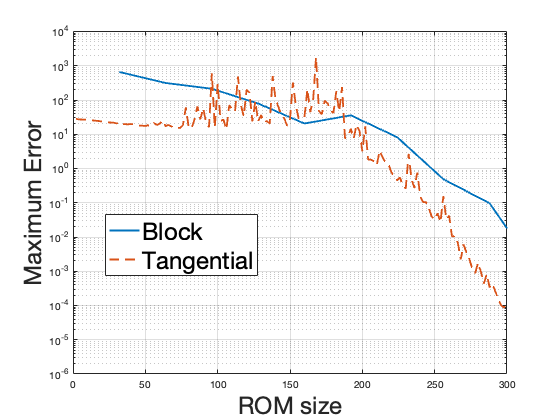}
\end{tabular}
\caption{{\bf Left:} Maximum Error Plots for Tangential and Block Adaptive Methods with 2 Sources per Side \\
{\bf Right:}  Maximum Error Plots for Tangential and Block Adaptive Methods with 4 Sources per side
}
\label{fig: max_error}
\end{figure}

We compared our approach with block adaptive algorithm that is similar to block-AAA approach. In block method one adds a full block $\bI_{m}$ of $m$ tangential vectors in one shot at step 7 of algorithm \ref{alg:adpol} rather than choosing them one by one in a greedy way. That allows Block method to achieve interpolation of the full transfer function, however, perhaps not all tangential vectors are needed if the goal is overall ROM accuracy. That said, when applying the Block and Tangential Adaptive algorithms to this medium, the two approaches differ in the number of basis vectors added to the reduced-order model at each selected frequency. For $m$ sources, the Block Adaptive method incorporates $2m$ vectors corresponding to each selected frequency, whereas the Tangential Adaptive method adds only two vectors, one associated with the negative frequency and one with the corresponding positive frequency. For the cases when not all tangential vectors equally contribute to the ROM accuracy,  the Tangential Adaptive method allows to achieve a superior accuracy for the same size of subspace. In particular, in our example the accuracy $10^{-4}$ is attained by our tangential approach with a smaller ROM size, while the Block Adaptive method performed significantly worse. These results indicate that the Tangential Adaptive method provides a more efficient approximation for this particular medium as seen in  Figure \ref{fig: max_error}.

\begin{figure}[t]
\centering
\hspace*{-1cm}
\renewcommand{\arraystretch}{1} 
\begin{tabular}{m{0.5cm} c @{\hspace{-0.4cm}} c} 
   & \includegraphics[height=4.67cm]{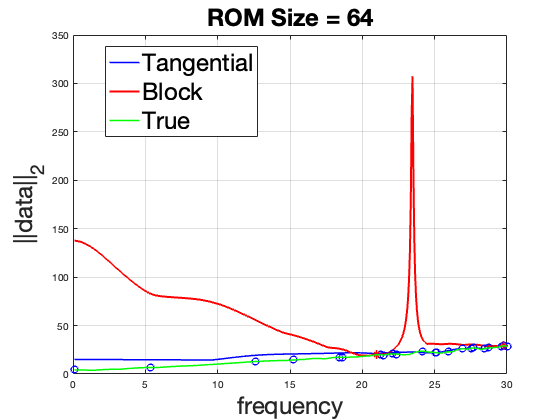} 
   & \includegraphics[height=4.65cm]{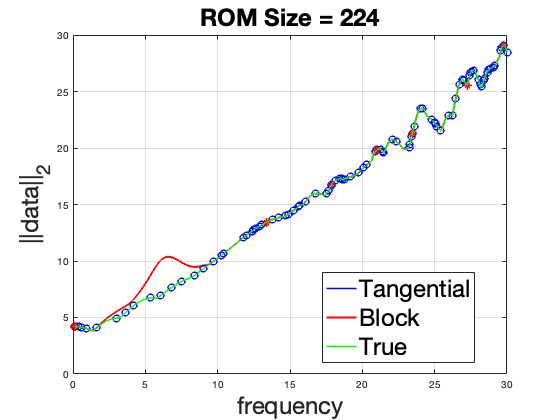} \\[0.2 cm]
   & \includegraphics[height=4.65cm]{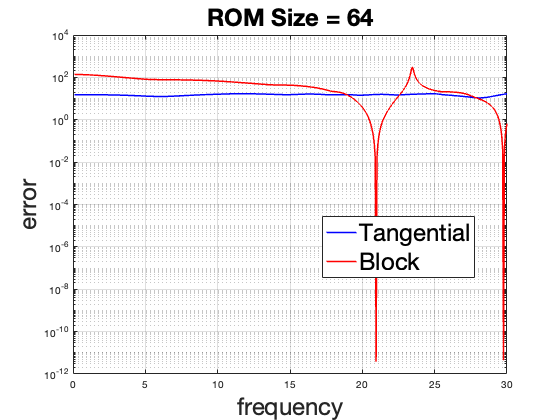} 
   & \includegraphics[height=4.65cm]{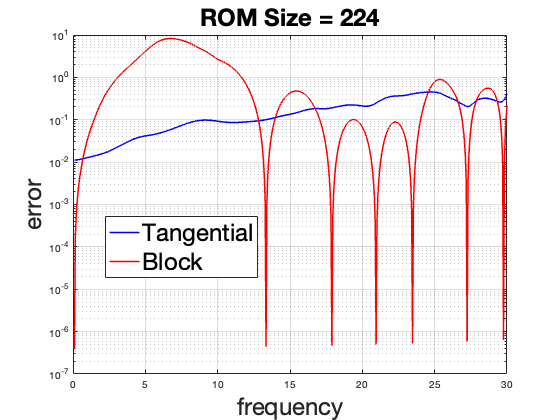} \\[0.2cm]
\end{tabular}
\caption{Comparison of the true and ROM transfer functions for ROM sizes 64 and 224. In the top two plots, the green curve represents the true transfer function, while the red and blue curves represent the ROM transfer functions obtained using the Block and Tangential Adaptive methods, respectively. Red stars and blue circles indicate the interpolation points selected by the Block and Tangential Adaptive methods, respectively. The bottom two plots show the corresponding errors between the true transfer function and the two ROM transfer functions.}

\label{fig: transferFunc}
\end{figure}

Comparing the transfer functions obtained using the two Adaptive algorithms using four sources per side, as shown in Figure \ref{fig: transferFunc}, we observe that the ROM transfer function generated by the Block Adaptive method interpolates the true transfer function only at  interpolation points. In contrast, the ROM transfer function produced by the Tangential Adaptive method provides a closer approximation to the true transfer function over the considered frequency range with a more uniform approximation error. This improved agreement results from the fact that, for the same ROM size, the Tangential Adaptive method selects more interpolation frequencies because it adds fewer basis vectors at each selected frequency.

\subsection{Spectral mass gauging for propellant tank monitoring} Spectral mass gauging (SMG) has been shown to be a promising tool to assess the volume of the remaining propellant in the tank of a  spaceship. It employs acoustic sources and receivers to assess the interior of the tank. Propagation of acoustic waves in the tank can be modelled in the frequency-domain via solving (\ref{dampedhelm}) from Example \ref{ex1} where $f^j=\delta(x-x^j)$ with $x^j$ being locations of collocated acoustic sources/receivers. Propellant tanks typically exhibit lossy behavior on the boundary only due to the foam used for insulation purposes. The volume of the propellant inside the tank in lossless scenario can be estimated via counting of resonances using Weyl's law. For 3D scenario can be formulated as $$vol(\Omega)=6\pi^2\lim_{\omega\rightarrow\infty}{\frac{N(\omega)}{\omega^3}}$$ where $N(\omega)$ is the number of resonances of (\ref{dampedhelm}) that are not greater than $\omega$. Although we consider propellant tank being lossy in our numerical experiment, damping is typically weak, so we can still adopt a similar volume approximation and count resonances with absolute value not greater than $\omega$. To assess the performance of our approach, in our second numerical example we consider the synthetic tank model $[0;\sqrt[4]{3}]\times [0;\sqrt[4]{2}]\times [0;1]$ that is known to be full. We generated synthetic data for 2000 frequencies from the interval $[0.9\pi;5.1\pi]$ for 6 collocated sources/receivers located on the boundary (one per face). To perform that, we employed numerical simulations using finite-difference method on $51\times 64\times 80$ 3D grid. Then we constructed data-driven ROM using our Algorithm \ref{alg:adpol} and assessed its performance to  count resonances via comparing them with the estimate given by Weyl's law. 

\begin{figure}[t]
\centering
\begin{tabular}{cc}
Transfer function & Resonances \\
\includegraphics[height=\figheight, keepaspectratio]{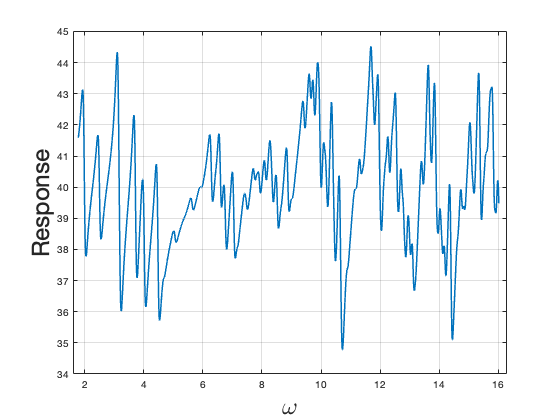} &
\includegraphics[height=\figheight, keepaspectratio]{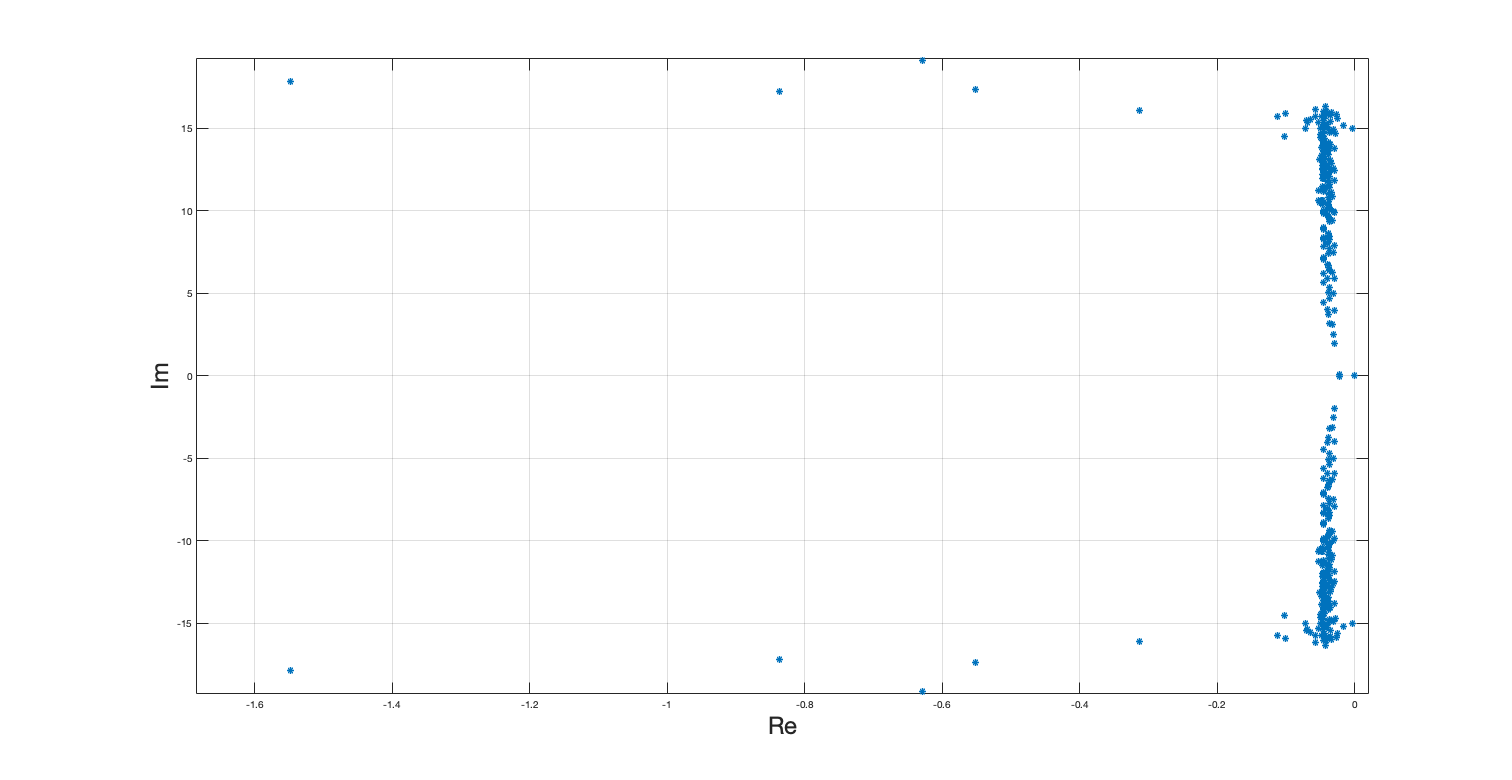}
\end{tabular}
\vspace{-1ex}
\caption{
{\bf Left:} Transfer function due to collocated source/receiver. 
{\bf Right:} Recovered resonances. Almost all of them have constant shift away from the imaginary axis. 
}

\label{resp&res}
\vspace{-2ex}
\end{figure}

In Figure \ref{resp&res} we plotted the data due to one of the sources recorded at the same location (left) and reconstructed resonances (right). As one can observe, most of the resonances have approximately the same shift from imaginary axis. In our experiments we
picked such boundary losses that shift is about 5\% of minimum positive imaginary part
of resonances.

\begin{figure}[t]
\centering
\begin{tabular}{c}
Counting function\\
\includegraphics[height=\figheight, keepaspectratio]{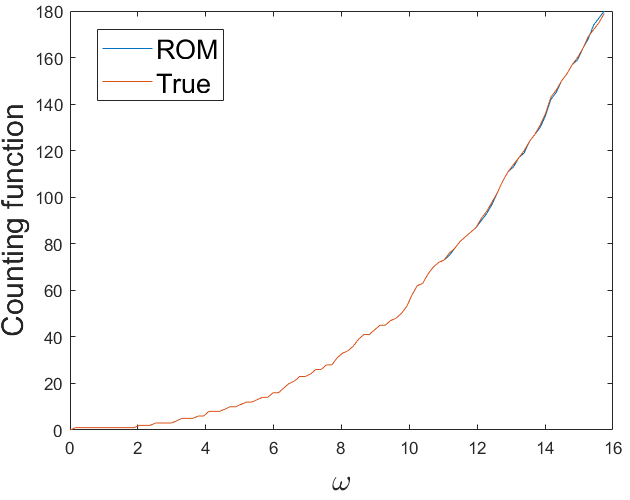}
\end{tabular}
\vspace{-1ex}
\caption{
Number of resonances counted by our ROM and counting function $N(\omega)$ obtained from Weyl's law.}

\label{cf}
\vspace{-2ex}
\end{figure}

In Figure \ref{cf} on the left we showed how our ROM managed to reproduce the
true counting function $N(\omega)$ for this model. In this example we successfully recovered about 180 resonances. That, in turn, allows us to infer the propellant volume with accuracy less than 1.5\% using Weyl's law.

\section{Conclusion}
\label{sec:conclusions}
Port-Hamiltonian systems provide a representation of a broad range of physical dynamical systems, but exploiting this structure in a fully data-driven reduced-order model is difficult when the underlying state-space realization is unavailable. In this work, we addressed this problem for the important subclass of reciprocal MIMO pH systems. Reciprocity supplies an additional symmetry that can be exploited directly at the reduced-order level: we showed that every minimal reciprocal realization is self-adjoint with respect to a generally indefinite inner product and that, up to a similarity transformation, this self-adjoint representation is compatible with the pH structure.

This hidden Krein-space structure provides the foundation for the proposed reduction procedure. A Galerkin projection with respect to the corresponding indefinite inner product produces a tangential rational Krylov ROM with a symmetric reduced matrix pencil and tangential Hermite interpolation properties. Although this model-driven construction formally requires full-order state solutions and knowledge of the indefinite inner product, we showed that both can be eliminated. The reduced mass and stiffness matrices can instead be constructed directly from samples of the transfer function and its frequency derivative through a Loewner-type formulation.

We further developed an adaptive data-driven procedure that jointly selects interpolation frequencies and tangential directions. Rather than interpolating the entire MIMO transfer function at every selected frequency, the method adds a direction associated with the dominant residual at the frequency where the current approximation error is largest. This allows the reduced basis to concentrate on those frequency regions and input-output directions that are most difficult to approximate. For the numerical examples considered here, the resulting tangential construction achieves substantially smaller approximation errors for a prescribed ROM dimension than block adaptive interpolation.

The spectral mass-gauging example further demonstrates that the resulting rational model can be used not only to approximate frequency-response data but also to extract spectral information. For the three-dimensional tank model, the data-driven ROM recovers approximately 180 resonances and produces a propellant-volume estimate with an error below 1.5\%  through Weyl's law. One possible further extension of applications range could be inverse scattering in dispersive media via data-driven ROMs \cite{zimmerling}. In particular, it would allow to adaptively pick the data samples that are the most influential to approximate the measured transfer function.

The present construction exploits reciprocity rather than attempting to solve the more general problem of data-driven structure-preserving reduction for arbitrary pH systems. At the same time, the resulting ROM does not in general preserve all properties of the original pH realization, in particular passivity and stability. Understanding how the reciprocal/Krein formulation developed here can be combined with additional constraints enforcing these properties is an important direction for future work. Other natural extensions include noisy or incomplete frequency-response data and adaptive experimental settings in which new frequency samples are acquired sequentially together with the construction of the ROM.

\thanks{{\bf Acknowledgments.} 
M. Zaslavsky was partially supported by  AFOSR grant  FA9550-23-1-0220. Authors are grateful to Vladimir Druskin, Mikhail Khasin, Cyrill Muratov and Joern Zimmerling for fruitful discussions that helped to improve this manuscript.}

\bibliographystyle{plain}
\bibliography{ROMbibliography}

\end{document}